\documentclass[12pt,reqno]{amsart}
\usepackage{amssymb,latexsym,amsmath,amsthm,mathrsfs,euscript,upgreek}
\usepackage{fullpage,enumerate}
\usepackage{amscd}
\usepackage[all]{xy}

\newcommand{\PP}{\mathbb P}
\newcommand{\QQ}{\mathbb Q}
\newcommand{\RR}{\mathbb R}
\newcommand{\ZZ}{\mathbb Z}

\newcommand{\al}{\alpha}
\newcommand{\gam}{\gamma}
\newcommand{\del}{\delta}
\newcommand{\lam}{\lambda}
\newcommand{\sig}{\sigma}

\newcommand{\Gam}{{\it \Gamma}}
\newcommand{\PHI}{{\it \Phi}}
\newcommand{\Sig}{{\it \Sigma}}

\newcommand{\XI}{{\it \Xi}}
\newcommand{\Th}{{\it \Theta}}
\newcommand{\Om}{{\it \Omega}}

\newcommand\g{\mathfrak{g}}
\newcommand\gl{\mathfrak{gl}}
\newcommand\su{\mathfrak{su}}

\def\cA{\EuScript A}

\def\cF{\EuScript F}
\def\cG{\EuScript G}
\def\cH{\EuScript H}

\def\cK{\EuScript K}
\def\cL{\EuScript L}
\def\cM{\EuScript M}
\def\cP{\EuScript P}
\def\cU{\EuScript U}

\def\bF{\mathbf{F}}
\def\bH{\mathbf{H}}

\def\bL{\mathbf{L}}

\def\Ad{\textup{Ad}}
\def\ad{\textup{ad}}
\def\id{\textup{id}}
\def\tr{\textup{tr}}
\def\Hom{\textup{Hom}}
\def\ker{\textup{ker}}
\def\coker{\textup{coker}}
\newcommand\irr{^\mathrm{irr}}
\def\bra{\langle}
\def\ket{\rangle}
\def\dd{\partial}
\def\ii{\sqrt{-1}}
\newcommand\set[2]{\{#1\,|\,#2\}}

\newcommand{\beq}{\begin{eqnarray}}
\newcommand{\beqn}{\begin{eqnarray*}}
\newcommand{\eeq}{\end{eqnarray}}
\newcommand{\eeqn}{\end{eqnarray*}}
\theoremstyle{plain}
\newtheorem{thm}{Theorem}[section]

\newtheorem*{thm*}{Theorem}

\newtheorem*{cor*}{Corollary}
\newtheorem*{conj*}{Conjecture}
\newtheorem*{lemma*}{Lemma}

\newtheorem*{prop*}{Proposition}

\newtheorem{proposition}[thm]{Proposition}
\theoremstyle{definition}

\newtheorem*{defn*}{Definition}

\newtheorem{examples*}{Examples}

\numberwithin{equation}{section}

\begin{document}

\title{Bundle gerbes and moduli spaces}

\author{Peter Bouwknegt}

\address[P Bouwknegt]{Department of Theoretical Physics,
Research School of Physics and
Engineering and Department of Mathematics, Mathematical Sciences Institute,
The Australian National University, Canberra, ACT 0200, Australia}
\email{peter.bouwknegt@anu.edu.au}

\author{Varghese Mathai}
\address[V Mathai]{Department of Pure Mathematics, University of Adelaide,
Adelaide, SA 5005, Australia}
\email{mathai.varghese@adelaide.edu.au}

\author{Siye Wu}
\address[S Wu]{Department of Mathematics, University of Hong Kong,
Pokfulam Road, Hong Kong, China}
\email{swu@maths.hku.hk}

\begin{abstract}
In this paper, we construct the {\em index bundle gerbe} of a family of
self-adjoint Dirac-type operators, refining a construction of Segal.
In a special case, we construct a geometric bundle gerbe called the 
{\em caloron bundle gerbe}, which comes with a natural connection and curving,
and show that it is isomorphic to the analytically constructed index bundle
gerbe.
We apply these constructions to certain moduli spaces associated to compact
Riemann surfaces, constructing on these moduli spaces, natural bundle gerbes
with connection and curving, whose 3-curvature represent Dixmier-Douady classes
that are generators of the third de Rham cohomology groups of these moduli
spaces.
\end{abstract}

\thanks{{\em Acknowledgements.}
This research was supported under Australian Research Council's Discovery
Projects funding scheme (project number DP0878184).
V.M.\ is the recipient of an Australian Research Council Australian
Professorial Fellowship (project number DP0770927).
S.W.\ is supported in part by CERG HKU705407P}


\subjclass[2000]{}
\maketitle

\section{Introduction}

Given a compact, simply-connected simple Lie group $G$, let $\cM$ denote the
moduli space of flat principal $G$-bundles on a compact Riemann surface $\Sig$
of genus greater than one.
More precisely, we only consider the smooth, dense, open subset of $\cM$
corresponding to irreducible homomorphisms from the fundamental group of
$\Sig$ to $G$, which we will denote by the same symbol.
Quantization of $\cM$ was considered by Witten~\cite{W,APW} via the
{\em determinant line bundle} $\cL$ of the index bundle of the associated
family of Cauchy-Riemann operators
$\{\bar\dd_\rho:\rho\in\Hom(\pi_1(\Sig),G)\}$.
The fibre of $\cL$ at $\rho$ is
$$  \cL_\rho=\bigwedge^{\rm max}\ker(\bar\dd_\rho^*)\otimes
    \bigwedge^{\rm max}\coker(\bar\dd_\rho).  $$
It carries the Quillen metric \cite{Q} and the curvature of the canonical
hermitian connection coincides with the natural K\"ahler form on the moduli
space $\cM$.
Moreover, the first Chern class of $\cL$ generates the second de Rham
cohomology group of $\cM$.
The moduli space $\cM$ plays an important role in quantum field theory.
In particular, the space of holomorphic sections of $\cL$ can be naturally
identified with what physicists call the space of {\em conformal blocks}
in a standard conformal field theory, the WZW-model~\cite{BL}, on which
there exists an extensive literature.

For a family of self-adjoint Dirac-type operators, an {\em index gerbe} was
introduced by Carey-Mickelsson-Murray~\cite{CM,CMM} and by Lott~\cite{L}.
The gerbe was described by local data, which can be viewed as the analogues
of the transition functions of the determinant line bundle.
The index gerbe is related to Hamiltonian anomalies.
In~\cite{S}, Segal constructs a projective Hilbert bundle; its Dixmier-Douady
invariant is the obstruction to lifting it to a Hilbert bundle.
In this paper, we construct the {\em index bundle gerbe} for such a family,
which is a global version of prior constructions.
Our construction of the index bundle gerbe is inspired by these papers and
by the work on determinant line bundles~\cite{BF,F,Q}.
Melrose-Rochon~\cite{MR} have a completely different construction of an index
bundle gerbe using pseudodifferential operators.
In the special case of a circle fibration, we also construct the
{\em caloron bundle gerbe} using geometric data and show that it
is isomorphic to the analytically defined index bundle gerbe.
We apply these constructions to the moduli space $\cM$ and obtain on it
natural bundle gerbes with connections and curvings.
The $3$-curvatures represent Dixmier-Douady classes that are generators
of the third de Rham cohomology group of $\cM$.

Applications of moduli spaces of Riemann surfaces are ubiquitous.
The construction of explicit geometric realizations of the degree $3$ classes
on moduli spaces through bundle gerbes was motivated by the desire to
generalize the geometric Langlands correspondence~\cite{Frenkel,HT,Hitchin}
by including background fluxes.

Here we outline the contents of the paper.
In Section~\ref{sect:index}, we construct the {\em index bundle gerbe}
associated to a set of geometric data using the spectrum of a family of
Dirac-type operators.
In Section~\ref{sect:caloron}, given a unitary representation of the
structure group, we construct the {\em caloron bundle gerbe} using the
projective Hilbert space bundle obtained from the caloron correspondence
of Murray-Stevenson~\cite{MS}.
In Section~\ref{sect:comparison}, we establish an isomorphism (not just a
stable isomorphism) between the geometrically constructed caloron bundle gerbe
and the analytically defined index bundle gerbe for a circle fibration.
Section~\ref{sect:moduli} contains the application of these constructions to
the moduli space $\cM$, obtaining natural bundle gerbes with connections and
curvings, whose Dixmier-Douady classes generate the third de Rham
cohomology group of $\cM$.
We end with a section of conclusions, where we outline future work on
constructing bundle gerbes on other moduli spaces such as the moduli space of
anti-self-dual connections, on a compact four dimensional Riemannian manifold.

\section{Index bundle gerbe}\label{sect:index}

In this section we will construct an index bundle gerbe associated to the
following set of geometric data.
For an introduction to bundle gerbes, see~\cite{Mur}.
\medskip

\noindent{\bf Basic setup}:\\
{\em Let $Z\to X$ be a smooth fibre bundle whose typical fibre is a compact
odd-dimensional manifold.
Let $T(Z/X)\to Z$ denote the vertical tangent bundle, which is a sub-bundle
of the tangent bundle $TZ$.
We assume that there is a Riemannian metric on $T(Z/X)$.
Suppose that $T(Z/X)$ has a spin structure and let $S\to Z$ denote the
corresponding bundle of spinors.
Let $E\to Z$ be a hermitian vector bundle with connection $\nabla^E$.}
\medskip

In this case, there is a smooth family $\{D^E_x:x\in X\}$ of self-adjoint
Dirac-type operators $D^E_x\colon\Gam(Z_x,S\otimes E)\to\Gam(Z_x,S\otimes E)$
acting on sections over the fibres $Z_x$ ($x\in X$) of $Z\to X$.
The $L^2$-completions $\cH_x$ of $\Gam(Z_x,S\otimes E)$ form a Hilbert
bundle $\cH$ over $X$.
We denote by the same symbol $D^E_x$ the operator on $\cH_x$.
Recall that for any $x\in X$, ${\rm spec}(D^E_x)\subset\RR$ is a closed
countable set with accumulation point only at infinity.
For $\lam\in\QQ$,\footnote{We restrict to $\lam\in\QQ$ to have a countable
open cover of $X$.} consider the open subset
\begin{equation}\label{eqn:eqBaa}
U_\lam=\{x\in X\big|\lam\notin{\rm spec}(D^E_x)\}
\end{equation}
of $X$.
For $x\in U_\lam$, let the Hilbert space $\cH_{\lam,x}^+$ ($\cH_{\lam,x}^-$,
respectively) be the span of eigenspaces of $D^E_x$ with eigenvalues greater
than (less than, respectively) $\lam$.
They form Hilbert bundles $\cH_{\lam}^\pm$ over $U_\lam$.

Suppose $\mu>\lam$ and $\lam,\mu\in\QQ$.
Observe that for $x\in U_\lam\cap U_\mu$, we have
$$ \cH_{\lam,x}^+=\cH_{\mu,x}^+\oplus(\cH_{\lam,x}^+\cap\cH_{\mu,x}^-),\quad
\cH_{\mu,x}^-=\cH_{\lam,x}^-\oplus(\cH_{\lam,x}^+\cap\cH_{\mu,x}^-), $$
where $\cH_{\lam,x}^+\cap\cH_{\mu,x}^-$ is a finite dimensional space.
Let $\cL_{\lam\mu}=\det(\cH_\lam^+\cap\cH_\mu^-)$, which is a line bundle
over $U_\lam\cap U_\mu$.
Since
$$   (\cH_\lam^+\cap\cH_\mu^-)\oplus(\cH_\mu^+\cap\cH_\tau^-)
=\cH_\lam^+\cap\cH_\tau^-   $$
on $U_\lam\cap U_\mu\cap U_\mu$ for any $\lam<\mu<\tau$ in $\QQ$, we have
a cocycle condition
\begin{equation}\label{cocycle}
\cL_{\lam\mu}\otimes\cL_{\mu\tau}\cong\cL_{\lam\tau}.
\end{equation}
The collection $\{\cL_{\lam\mu}\}$ defines the index gerbe~\cite{CM,CMM,L}
over $X$.

Over $U_\lam$, we construct the Fock bundle
$\cF_\lam=\bigwedge(\cH_{\lam}^{+}\oplus\overline{\cH_{\lam}^{-}})$.
Hodge duality asserts that for any finite dimensional complex vector space $V$
there is a canonical isomorphism
$$ \bigwedge\overline V\otimes\det V\cong\bigwedge V. $$
Thus on the overlap $U_\lam\cap U_\mu$, we get
\begin{align*}\label{eqn:eqBa}
\cF_\lam
&=\bigwedge\cH_\lam^+\otimes\bigwedge\overline{\cH_\lam^-}     \\
&\cong\bigwedge\cH_\mu^+\otimes\bigwedge(\cH_\lam^+\cap\cH_\mu^-)
 \otimes\bigwedge\overline{\cH_\lam^-}                         \\
&\cong\bigwedge\cH_\mu^+\otimes\bigwedge\overline{\cH_\lam^-}\otimes
 \bigwedge\overline{(\cH_\lam^+\cap\cH_\mu^-)}\otimes\cL_{\lam\mu} \\
&\cong\bigwedge\cH_\mu^+\otimes\bigwedge\overline{\cH_\mu^-}
 \otimes\cL_{\lam\mu}                                          \\
&=\cF_\mu\otimes\cL_{\lam\mu}.
\end{align*}
Therefore there is a canonical isomorphism of the projectivizations
$\PP(\cF_\lam)\cong\PP(\cF_\mu)$ on the overlap $U_\lam\cap U_\mu$ and
hence a well defined projectivized Fock bundle $\pi\colon\PP(\cF)\to X$.

To construct the index bundle gerbe, we now define a line bundle
$\cL\to\PP(\cF)^{[2]}$ on the fibre product of $\PP(\cF)$ with itself.
Let $\cL_\lam$ be the universal line bundle over
$\PP(\cF_\lam)=\pi^{-1}(U_\lam)$.
That is, over any $\ell\in\PP(\cF_\lam)$, the fibre is
$$ (\cL_\lam)_\ell=\{(\ell,v):v\in\ell\}\subset
   \{\ell\}\times(\cF_\lam)_{\pi(\ell)}. $$
$\cL_\lam$ is a hermitian line bundle, but does not glue properly on the
overlaps in $\PP(\cF)$ for the same reason that $\cF_\lam$ doesn't on $X$.
In fact, on $\pi^{-1}(U_\lam)\cap\pi^{-1}(U_\mu)=\pi^{-1}(U_\lam\cap U_\mu)$,
we have
\begin{equation}\label{cocycle2}
\cL_\lam \cong \cL_\mu\otimes\pi^*\cL_{\lam\mu}.
\end{equation}
Now, $\del(\cL_\lam)=\cL_\lam^*\boxtimes\cL_\lam$ is a line bundle on
$\PP(\cF_\lam)\times\PP(\cF_\lam)$ which restricts to the total space of
the bundle $\pi_\lam^{[2]}\colon\PP(\cF_\lam)^{[2]}\to U_\lam$.
However, this time, $\del(\cL_\lam)=\del(\cL_\mu)$ on the overlap
$(\pi_\lam^{[2]})^{-1}(U_\lam\cap U_\mu)$, since the fudge factor for each
cancels, thereby defining a global hermitian line bundle
$\cL\to\PP(\cF)^{[2]}$.
There is a line bundle $\del(\cL)\to\PP(\cF)^{[3]}$ defined by
$$ \del(\cL)=\pi_{12}^*\cL\otimes(\pi_{13}^*\cL)^{-1}\otimes\pi_{23}^*\cL, $$
where $\pi_{ij}\colon\PP(\cF)^{[3]}\to\PP(\cF)^{[2]}$ for $1\le i<j\le 3$ is
the projection onto the $i$th and $j$th factors.
It is trivial precisely because of~\eqref{cocycle}.

Therefore we have established the following.

\begin{thm}
Given the {\bf Basic setup}, there is a well-defined global fibre bundle
of projective Hilbert spaces $\pi\colon\PP(\cF)\to X$ and a line bundle
$\cL\to\PP(\cF)^{[2]}$ such that $\del(\cL)\to\PP(\cF)^{[3]}$ is trivial.
\end{thm}

Thus we have a bundle gerbe over $X$, which we call the {\em index bundle
gerbe}.
Note that the line bundle $\cL_\lam\to\PP(\cF_\lam)$ provides a local
trivialization of the index bundle gerbe over $U_\lam\subset X$.

We remark that the considerations in this section can straightforwardly be
generalized to more general open coverings $\{U_\al\}$ of $X$ so that for each
$U_\al$ there exists a spectral cut $\lam_\al\notin\text{spec}(D_x^E)$, for
all $x\in U_\al$.
We restrict our discussion to the cover~\eqref{eqn:eqBaa} for notational
simplicity.
Finally, we remark that the construction in this section works for a family
of general elliptic self-adjoint operators on a compact manifold.

\section{The caloron bundle gerbe}\label{sect:caloron}

We begin by reviewing the caloron correspondence of Murray-Stevenson \cite{MS}.
We then study the behavior of this construction under a unitary representation
of the structure group.
Finally, we obtain a fibre bundle whose typical fibre is a projective Hilbert
space; we will call the associated bundle gerbe the caloron bundle gerbe.

\subsection{The caloron correspondence}\label{sect:calorona}

Let $X$ be a manifold and let $S^1$ denote the unit circle.
Let $P\to S^1\times X$ be a principal $G$-bundle, where $G$ is a connected
and simply-connected compact Lie group.
Then there exists a canonical identification between equivalence classes of
principal $G$-bundles over $S^1\times X$ and principal $LG$-bundles over $X$,
\begin{equation} \label{eqBa}
\begin{pmatrix} & \begin{CD}
G @>>> P \\ && @VV V \\ && S^1 \times X
\end{CD} & \end{pmatrix}
\quad\Leftrightarrow\quad
\begin{pmatrix} & \begin{CD}
LG @>>> Q \\  && @VV V \\   && X
\end{CD} & \end{pmatrix},
\end{equation}
as follows.
Given $P$, the bundle $Q$ is the push-forward of $P$ under the projection
map $S^1\times X\to X$.
That is, the fibre of $Q$ over $x\in X$ is $Q_x=\Gam(S^1\times\{x\},P)$.
Clearly, the free loop group $LG$ acts freely and transitively on $Q_x$ and
hence $Q$ is a principal $LG$-bundle.
The total space $Q$ can be constructed globally as follows.
Consider the $LG$-bundle $LP\to L(S^1\times X)$.
Then $Q\to X$ is the pullback of $LP$ by the map
$\eta\colon X\to L(S^1\times X)$ given by $\eta(x)=(\theta\mapsto(\theta,x))$.
Conversely, given an $LG$-bundle $Q\to X$, we can define a principal
$G$-bundle $P$ over $S^1\times X$ by $P=(Q\times G\times S^1)/LG$, where the
right $LG$-action on $Q\times G\times S^1$ is given by
$(p,g,\theta)\gam=(p\gam,\gam(\theta)^{-1}g,\theta)$
and the $G$-action is by right multiplication in the second factor.
Clearly these actions commute, establishing the canonical equivalence in
\eqref{eqBa}.
For a generalization of the caloron correspondence to principal $G$-bundles
over non-trivial circle bundles or over more general fibrations,
see~\cite{BV,BM,MV,HMV}.

Given a connection $\tilde{A}$ on $P$, we define a connection
$A=\eta^*L\tilde A$ on $Q$ by pulling back the connection $L\tilde A$ on $LP$.
In addition, $\tilde A$ determines a Higgs field
$\PHI=\eta^*\iota_\XI(L\tilde A)$, where $\XI$ is a canonical vector field
on $LP$ generating the translation on $S^1$.
The map $\PHI\colon Q\to L{\mathfrak g}$ is smooth and satisfies
\begin{equation} \label{eqBb}
\PHI(p\gam)=\ad(\gam^{-1})\PHI(p)+\gam^{-1}\dd_\theta\gam  
\end{equation}
for $\gam\in LG$.
Conversely, given a connection $A$ on $Q$ and a Higgs field $\PHI$ for $Q$
satisfying Eqn.~\eqref{eqBb}, we have a $1$-form
$$ \tilde A=\ad(g^{-1})A(\theta)+\Th+\ad(g^{-1})\PHI\,d\theta $$
on $Q\times G\times S^1$ which descends to a connection $1$-form on $P$.
Here $\Th$ denotes the Cartan-Maurer $1$-form on $G$.
We summerize the above in the following

\begin{proposition}[caloron correspondence \cite{MS}]
The canonical equivalence in \eqref{eqBa} determines a bijection between
isomorphism classes of $G$-bundles with a connection over $S^1\times X$ and
isomorphism classes of $LG$-bundles with a connection and a Higgs field
over $X$.
\end{proposition}

Recall that for a simple group $G$, there is a basic central extension
(cf.~\cite{PS})
$$ 1\longrightarrow U(1)\longrightarrow\widehat{LG}
\stackrel{p}{\longrightarrow}LG\longrightarrow 1 $$
corresponding to a cocycle that generates the second cohomology of $LG$.
Denote by $\xi\colon Q^{[2]}\to LG$ the map given by
$(q_1,q_2)\mapsto\gam$, where $\gam\in LG$ is the unique element such that
$q_2=q_1\gam$.
Let $\cP^Q\to Q^{[2]}$ be the pullback of the $U(1)$-bundle
$\widehat{LG}\to LG$ via $\xi$ and let $\cL^Q\to Q^{[2]}$ be
the line bundle associated to $\cP^Q$.
Equivalently, if $\cL^{LG}\to LG$ is the line bundle associated to
$\widehat{LG}$, then $\cL^Q=\xi^*(\cL^{LG})$.
Since $\xi$ is a homomorphism of groupoids, the line bundle $\cL^Q$ determines
a bundle gerbe called the {\it lifting bundle gerbe} of $Q$.

Theorem 5.1 of \cite{MS} asserts that a connection $A$ on $Q$ and a Higgs
field for $Q$ determine a curving (or B-field)
$$ B=\frac{\ii}{2\pi}\int_{S^1}\left(\frac{1}{2}\bra A,\dd_\theta A
   \ket-\bra F_A,\nabla\PHI\ket \right)d\theta,  $$
where $F_A$ is the curvature of the connection $A$ and
$\nabla\PHI=d\PHI+[A,\PHI]-\dd_\theta A$.
Here and below, the inner product $\bra\cdot,\cdot\ket$ on the Lie algebra
${\mathfrak g}$ is chosen such that the length of a long root is $\sqrt{2}$.
It induces an inner product (with the same notation) on $L{\mathfrak g}$.
The $3$-curvature $H=dB\in\Om^3(X)$, associated to $B$, represents the
Dixmier-Douady class of the bundle gerbe $(Q,\cL^Q)$ in de Rham cohomology.
It pulls back to
$-\frac{\ii}{4\pi^2}\int_{S^1}\bra F_A,\nabla\PHI\ket\,d\theta$ on $Q$.
Let $F_{\tilde A}$ denote the curvature of the connection $\tilde A$ on $P$.
Then the first Pontryagin form of $P$ is
$\frac{1}{8\pi^2}\bra F_{\tilde A},F_{\tilde A}\ket$, and Theorem 6.1 
in~\cite{MS} asserts that
\begin{equation}\label{MSthm}
-\frac{1}{8\pi^2}\int_{S^1}\bra F_{\tilde A},F_{\tilde A}\ket=H.
\end{equation}

If $\rho\colon G\to SU(V)$ is a unitary representation of $G$ on a finite
dimensional vector space $V$, we can study the behavior of the caloron
correspondence under $\rho$.
First, we have a principal $SU(V)$-bundle
$P^\rho=P\times_GSU(V)\to S^1\times X$;
let $r_P\colon P\to P^\rho$ be the map induced by $\rho$ that changes
the structure group.
Correspondingly, there is a principal $LSU(V)$-bundle
$Q^\rho=\eta^*(LP^\rho)=Q\times_{LG}LSU(V)$ over $X$.
Let $r_Q\colon Q\to Q^\rho$ be the map induced by $L\rho\colon LG\to LSU(V)$.
Using the map $\xi_\rho\colon(Q^\rho){}^{[2]}\to LSU(V)$, we have a line
bundle $\cL^\rho=\xi_\rho^*(\cL^{LSU(V)})$ and hence a bundle gerbe
$(Q^\rho,\cL^\rho)$.
We recall that the line bundles $\cL^{LG}\to LG$ and $\cL^{LSU(V)}\to LSU(V)$
satisfy the relation
$(L\rho)^*(\cL^{LSU(V)})\cong(\cL^{LG})^{\otimes\iota_\rho}$, where
$\iota_\rho$ is the Dynkin index of the representation $\rho$ \cite{BETZ}.
Using the commutative diagram
$$ \xymatrix{Q^{[2]} \ar[r]^\xi \ar[d]_{r_Q^{[2]}} & LG \ar[d]^{L\rho} \\
(Q^\rho)^{[2]} \ar[r]^{\xi_\rho} & LSU(V),}    $$
we get $(r_Q^{[2]})^*(\cL^\rho)\cong(\cL^Q)^{\otimes\iota_\rho}$.
Therefore, the bundle gerbe $(Q^\rho,\cL^\rho)$ on $X$ is the $\iota_\rho$-th
power of $(Q,\cL^Q)$.

Given a connection $\tilde A$ of $P\to S^1\times X$, the induced connection
$\tilde A_\rho$ on $P^\rho$ satisfies
$r_P^*(\tilde A_\rho)=\dot\rho(\tilde A)$, where $\dot\rho\colon\g\to\su(V)$
is the representation of the Lie algebra $\g$.
The connection on $Q^\rho$ is $A_\rho=\eta^*(L\tilde A_\rho)$; it satisfies
$r_Q^*(A_\rho)=L\dot\rho(A)$, and therefore we have
$r_Q^*(F_{A_\rho})=L\dot\rho(F_A)$ for the corresponding curvatures.
Similarly, the Higgs field
$\PHI_\rho=\eta^*\iota_\XI(L\tilde A_\rho)\colon Q^\rho\to L\su(V)$ satisfies
$r_Q^*(\PHI_\rho)=L\dot\rho(\PHI)$.
Since the pull-back by $\dot\rho$ of the primitive bilinear form on $\su(V)$
is $\iota_\rho$ times the form $\bra\cdot,\cdot\ket$ on $\g$, the curving
(B-field) $B_\rho$ on $Q^\rho$ is related to $B$ on $Q$ by
$r_Q^*(B_\rho)=\iota_\rho\,B$.
Therefore, the $3$-curvature $H_\rho$ on $X$ associated to $B_\rho$ is
$H_\rho=\iota_\rho\,H$.
This is compatible with an earlier result that the bundle gerbe
$(Q^\rho,\cL^\rho)$ is the $\iota_\rho$-th power $(Q,\cL^Q)$.

We summarize the results in the following

\begin{proposition}\label{hom}
Given a representation $\rho\colon G\to SU(V)$ in the above setting, the
induced bundle gerbe $(Q^\rho,\cL^\rho)$ on $X$ is the $\iota_\rho$-th power
of $(Q,\cL^Q)$. 
Furthermore, its B-field and the $3$-curvature are, respectively,
$$ r_Q^*B_\rho=\iota_\rho\,B,\qquad H_\rho=\iota_\rho\,H. $$
\end{proposition}

\subsection{Projective Hilbert bundle and the caloron bundle gerbe}
\label{sect:caloronb}

As before, $G$ is a simple, connected and simply-connected compact Lie group
and $\rho\colon G\to SU(V)$ is a unitary representation of $G$.
Let $\bH=L^2(S^1,V)$.
Then there is an induced unitary representation $\tilde\rho\colon LG\to U(\bH)$
of $LG$ defined by
$\tilde\rho(\gam)(\theta)=\rho(\gam(\theta))$, where $\gam\in LG$ and
$\theta\in S^1$.
For any $\lam\in\RR$, we have a decomposition
$\bH=\bH_\lam^+\oplus\bH_\lam^-$, where $\bH_\lam^+$ ($\bH_\lam^-$,
respectively) denotes the Hilbert space span of the Fourier modes not less than
(less than, respectively) $\lam$.
Actually, we have $\tilde\rho\colon LG\to U_{\rm res}(\bH)$, where the
{\em restricted unitary group} is
$$ U_{\rm res}(\bH)=\{S\in U(\bH):[S,I_\lam]\,\,
   \text{is a Hilbert-Schmidt operator}\}\,.   $$
Here $I_\lam$ is an involution on $\bH$ such that $I_\lam=\id$ on $\bH_\lam^+$
and $I_\lam=-\id$ on $\bH_\lam^-$ (cf.~\cite{PS}).
It is clear that $U_{\rm res}(\bH)$ does not depend on the choice of $\lam$.

Let $\bF_\lam$ denote the associated Fock (Hilbert) space
$$  \bF_\lam=\bigwedge(\bH_\lam^+\oplus\overline{\bH_\lam^-})\,. $$
Then we have the Shale-Stinespring embedding
$\sig_\lam\colon U_{\rm res}(\bH)\hookrightarrow PU(\bF_\lam)$
(cf.~\cite{Wurz}).
Let $\hat\rho_\lam=\sig_\lam\circ\tilde\rho\colon LG\to PU(\bF_\lam)$
be the composition.
It is a positive energy representation of $LG$.
In the next section, we will compare the representations for different values
of $\lam$, but for the remainder of this section we will only need $\lam=0$.

Given a principal $G$-bundle $P\to S^1\times X$, we have the principal
$LG$-bundle $Q\to X$ by the caloron correspondence.
Together with a representation $\rho$ of $G$ on $V$, we form the associated
fibre bundle
$$ \PP(\cF^\rho)=Q\times_{LG}\PP(\bF_0)\to X, $$
where $LG$ acts on the typical fibre $\PP(\bF_0)$ via $\hat\rho_0$.
We will next define a line bundle $\cL^\rho\to\PP(\cF^\rho)^{[2]}$ induced by
$\cP^Q$ or $\cL^Q\to Q^{[2]}$.
Since
$$ \cP^Q=\set{(q_1,q_2,g)\in Q^{[2]}\times\widehat{LG}}{q_2=q_1p(g)}, $$
there is a right action of $\widehat{LG}\times\widehat{LG}$ on $\cP^Q$ given by
$$ (g_1,g_2)\colon(q_1,q_2,g)\mapsto(q_1p(g_1),q_2p(g_2),g_1^{-1}gg_2). $$
On the other hand, if we denote the universal line bundle over $\PP(\bF_0)$
by $\bL$, the left $LG$-action on $\PP(\bF_0)$ (via $\hat\rho_0$) lifts to
an action of $\widehat{LG}$ on $\bL$.
We set
$\cL^\rho=\cP^Q\times_{\widehat{LG}\times\widehat{LG}}(\bL\boxtimes\bL)$;
this is a line bundle over $Q^{[2]}\times_{(LG\times LG)}
(\PP(\bF_0)\times\PP(\bF_0))=\PP(\cF^\rho)^{[2]}$.
The triviality of $\del(\cL^Q)\to Q^{[3]}$ implies that of
$\del(\cL^\rho)\to\PP(\cF^\rho)^{[3]}$.
Thus we have a bundle gerbe which we call the {\em caloron bundle gerbe}.

We explain the above construction using local data.
Let $\{U_\al\}$ be an open cover of $X$ which locally trivializes the
principal $LG$-bundle $Q$.
Then the restriction of $Q$ to $U_\al$, $Q_\al\to U_\al$, has a lift
to a principal $\widehat{LG}$-bundle $\hat Q_\al\to U_\al$.
On the overlap $U_\al\cap U_\beta$, the two bundles $\hat Q_\al$ and
$\hat Q_\beta$ differ by a principal $U(1)$-bundle whose associated line
bundle we denote by $\cL^Q_{\al\beta}\to U_\al\cap U_\beta$.
They satisfy $\cL^Q_{\al\beta}\otimes\cL^Q_{\beta\gam}\cong\cL^Q_{\al\gam}$
on $U_\al\cap U_\beta\cap U_\gam$ and the collection $\{\cL^Q_{\al\beta}\}$
is a local description of the bundle gerbe $(Q,\cL^Q)$.
Given a unitary representation $\rho\colon G\to SU(V)$, the principal
$LSU(V)$-bundle $Q^\rho_\al=Q_\al\times_{LG}LSU(V)$ lifts to a principal
$\widehat{LSU(V)}$-bundle
$\widehat{Q^\rho_\al}=\hat Q_\al\times_{\widehat{LG}}\widehat{LSU(V)}$
over $U_\al$.
On the overlap $U_\al\cap U_\beta$, the two lifts $\widehat{Q^\rho_\al}$ and
$\widehat{Q^\rho_\beta}$ differ by a $U(1)$-bundle whose associated line
bundle is $\cL^\rho_{\al\beta}=(\cL^Q_{\al\beta})^{\otimes\iota_\rho}$.

Let $\cF^\rho_\al=\hat Q_\al\times_{\widehat{LG}}\bF_0\to U_\al$ be the
associated Fock bundles; the projectivizations $\PP(\cF^\rho_\al)$ glue
properly to form the bundle $\pi_\rho\colon\PP(\cF^\rho)\to X$.
On the overlap $U_\al\cap U_\beta$, the Fock bundles are related by
$\cF^\rho_\al\cong(\cL^Q_{\al\beta})^{\otimes\iota_\rho}\otimes\cF^\rho_\beta$.
As in the construction of the index bundle gerbe, let $\cL^\rho_\al$ be the
universal bundle over $\PP(\cF^\rho_\al)$.
Then $\cL^\rho_\al=\cL^\rho_\beta\otimes\pi_\rho^*(\cL^\rho_{\al\beta})$.
So the line bundles $(\cL^\rho_\al)^{-1}\boxtimes\cL^\rho_\al$ over
$\PP(\cF^\rho_\al)^{[2]}\subset\PP(\cF^\rho_\al)\times\PP(\cF^\rho_\al)$ glue
properly to define globally the line bundle $\cL^\rho\to\PP(\cF^\rho)^{[2]}$.
This defines the same caloron bundle gerbe. 
The line bundles $\cL^\rho_\al$ are its local trivializations under which the
bundle gerbe is described locally by a collection of line bundles
$\cL^\rho_{\al\beta}=(\cL^Q_{\al\beta})^{\otimes\iota_\rho}$ over
$U_\al\cap U_\beta$.

\section{Comparison of the index bundle gerbe and the caloron bundle gerbe}
\label{sect:comparison}

Recall that $P$ is a principal $G$-bundle over $S^1\times X$ and $Q$ is the
corresponding principal $LG$-bundle over $X$.
The fibre of $Q$ at $x\in X$ is $Q_x=\Gam(S^1\times\{x\},P)$ with the obvious
right $LG$-action.
Given a finite dimensional unitary representation $\rho$ of $G$ on $V$, there
is an associated hermitian vector bundle $E=P\times_GV\to S^1\times X$.
Consider the Hilbert space bundle $\cH\to X$ whose fibre $\cH_x$ over $x\in X$
is the $L^2$-completion of the space $\Gam(S^1\times\{x\},E)$.
There is a family of Dirac operators $\{D_x\}$ on $S^1$ coupled to $E$ acting
on the fibres of $\cH$.
By Sect.~\ref{sect:index}, there is a projectivized Fock bundle $\PP(\cF)\to X$
and a hermitian line bundle $\cL\to\PP(\cF)^{[2]}$ that defines the index
bundle gerbe.
By Sect.~\ref{sect:caloronb}, there is a another projectivized Fock bundle 
$\PP(\cF^\rho)\to X$ constructed by using the data from caloron correspondence
and the spectral cut at $0$ and a hermitian line bundle
$\cL^\rho\to\PP(\cF^\rho)^{[2]}$ that defines the caloron bundle gerbe.
The purpose of this section is to show that the two bundle gerbes are
isomorphic.

\subsection{Spectral cuts, equivalent representations and Bogoliubov
transformations}\label{sect:Bog}

In this section, we show that the projective loop group representation the
projectivized Fock space $\PP(\bF_\lam)$, corresponding to a spectral cut
$\lam$, is equivalent to $\PP(\bF_\mu)$, for $\mu\neq\lam$.
This is achieved by using a Bogoliubov transformation 
(see also \cite{Mick}, Sect.~12.3).

Let $\{\psi^i_n:n\in\ZZ,1\le i\le\dim V\}$ be a basis for $\bH$.
Then $\{\bar\psi^i_n:n\in\ZZ,1\le i\le\dim V\}$ is a basis for $\overline\bH$,
where $\bar\psi^i_n$ is the complex conjugation of $\psi^i_{-n}$.
They are chosen so that $\bH_\lam^+$ is spanned by
$\{\psi^i_n:n\in\ZZ,n\ge\lam,1\le i\le\dim V\}$ and $\bH_\lam^-$
is spanned by $\{\psi^i_n:n\in\ZZ,n<\lam,1\le i\le\dim V\}$.
Similarly, $\overline{\bH_\lam^+}$ is spanned by
$\{\bar\psi^i_n:n\in\ZZ,n\le-\lam,1\le i\le\dim V\}$ and
$\overline{\bH_\lam^-}$ is spanned by
$\{\bar\psi^i_n:n\in\ZZ,n>-\lam,1\le i\le\dim V\}$.
The Fock space $\bF_\lam$ is spanned by
$\{\psi^{i_1}_{n_1}\ldots\psi^{i_N}_{n_N}\bar\psi^{j_1}_{m_1}\ldots
\bar\psi^{j_M}_{m_M}|\lam\ket\ \,,\ n_1,\ldots,n_N\ge\lam,\
m_1,\ldots,m_M>-\lam\}$, where the `vacuum' $|\lam\ket$ has the properties
\begin{equation*}
\psi^i_n\,|\lam\ket=0\,,\quad\,n<\lam\,\quad\text{and}\quad
\bar\psi^i_n\,|\lam\ket=0\,,\quad \, n\le-\lam\,.
\end{equation*}
and $\psi_m^i$ and $\bar\psi_n^j$ are operators act on $\bF_\lam$ satisfying
\begin{align*}
\{\psi^i_m,\psi^j_n\}&=0=\{\bar\psi^i_m,\bar\psi^j_n\}\,,  \\
\{\psi^i_m,\bar\psi^j_n\}&=\del_{m+n,0}\del^{i,j}\,.
\end{align*}

If $\mu>\lam$, then
\begin{equation*}
|\mu\ket=\left({\prod_{\lam<n\le\mu}}\prod_i\psi^i_n\right)
|\lam\ket\,,\qquad
|\lam\ket=\left({\prod_{-\mu\le n<-\lam}}\prod_i\bar\psi^i_n\right)
|\mu\ket\,.
\end{equation*}
The operator $\displaystyle\prod_{\lam<n\le\mu}\prod_i\psi^i_n$, corresponding
to a base vector of the complex line $\det(\bH^{\lam+}\cap\bH^{\mu-})$, relates
the two vacua in $\bF_\lam$ and $\bF_\mu$; this is the so called Bogoliubov
transformation.

Just as $U_{\rm res}(\bH)$, we have the restricted general linear group
$$ GL_{\rm res}(\bH)=\{U\in GL(\bH):[U,I_\lam]\,\,
   \text{is a Hilbert-Schmidt operator}\}\,,   $$
which is also independent of $\lam$.
Define the bounded operator $e^{ij}_n$ on $\bH$ by
$e^{ij}_n(\psi^k_m)=\del^{j,k}\psi^i_{m+n}$.
Then $\id+e^{ij}_n\in GL_{\rm res}(\bH)$.
On the Fock space $\bF_\lam$ we define the operators
$\sig_\lam(e^{ij}_n)=\sum_m\!:\!\psi^i_m\bar\psi^j_{n-m}\!:_\lam\,$,
where the $\lam$-normal ordering is defined as
$$:\!\psi^i_m \bar\psi^j_{n}\!:_\lam=\left\{
\begin{array}{lcl}\psi^i_m\bar\psi^j_{n},\quad\mbox{if }m\le\lam\,, \\
-\bar\psi^j_{n}\psi^i_m,\quad\mbox{if }m<\lam\,. \end{array} \right.
$$
A standard computation then gives
\begin{equation*}
[\sig_\lam(e^{ij}_m),\sig_\lam(e^{kl}_n)]=
\del^{j,k}\sig_\lam(e^{il}_{m+n})-\del^{i,l}\sig_\lam(e^{kj}_{m+n})
+\del_{j,k} \del^{i,l} m \del_{m+n,0}\,,
\end{equation*}
which shows that $\sig_\lam$ is a representation on $\cF_\lam$ of the
central extension $\widehat{\gl_{\rm res}(\bH)}$ of the Lie algebra 
$\gl_{\rm res}(\bH)$ of $GL_{\rm res}(\bH)$.

Now let $\mu>\lam$.
Then one computes
\begin{align*}
\sig_\mu(e^{ij}_n)&=\sig_\lam(e^{ij}_n)-\sum_{\lam<m<\mu}
                                    \{\psi^i_m,\bar\psi^j_{n-m}\}       \\
&=\sig_\lam(e^{ij}_n)-n_{\lam\mu}\del^{i,j}\del_{n,0},
\end{align*}
where $n_{\lam\mu}=\#\{m:\lam<m\le\mu\}$.
Since $e^{ij}_n$ generates $\widehat{\gl_{\rm res}(\bH)}$, we have
$$ \sig_\mu(K)=\sig_\lam(K)-n_{\lam\mu}{\rm Tr}(K) $$
for any $K\in\widehat{\gl_{\rm res}(\bH)}$ and hence
$$ \sig_\mu(U)=\sig_\lam(U)\det(U)^{-n_{\lam\mu}} $$
for any $U\in\widehat{GL_{\rm res}(\bH)}$.
That is, the projective representations on $\PP(\bF_\lam)$ and $\PP(\bF_\mu)$
of $U_{\rm res}(\bH)\subset GL_{\rm res}(\bH)$, and hence those of $LG$
are equal.

We note that in this particular case, not only do we have an action of a
central extension of $GL_{\text{res}}(\bH)$ on $\bF_\lam$, but also an
action of the Virasoro algebra (the central extension of the 2-dimensional
conformal algebra).
The calculations in this section are of course well-known to the experts in
conformal field theory (see, e.g.,~\cite{DFMS} and references therein).
It is an interesting question whether our construction of the caloron bundle
gerbe has applications in the context of conformal field theory as well.

\subsection{Isomorphism between the index and caloron bundle gerbes}
\label{theisom}

In this section, we will show that the geometrically obtained caloron bundle
gerbe is isomorphic (not just stably isomorphic) to the analytically defined
index bundle gerbe in this context.
For the notion of stable isomorphism of bundle gerbes, see~\cite{MS2}.

Consider the vector bundle $E=P\times_GV$ over $S^1\times X$, which
itself fibers over $X$ with circle fibres.
We first show that Hilbert bundle $\cH$ defined in Sect.~\ref{sect:index} is
the associated bundle of $Q$ by $\tilde\rho$, i.e., $\cH\cong Q\times_{LG}\bH$.
Suppose $X$ has an open covering $\{U_\al\}$ such that the bundle $Q$ is
trivial over each $U_\al$.
Then $P$ is also trivial on $S^1\times U_\al$; let
$g_{\al\beta}\colon S^1\times(U_\al\cap U_\beta)\to G$ be the
transition functions.
Then the transition functions
$\tilde g_{\al\beta}\colon U_\al\cap U_\beta\to LG$ of $Q$ are given by
$\tilde g_{\al\beta}(x)=g_{\al\beta}(\cdot,x)$, $x\in U_\al\cap U_\beta$.
Under the trivialisation over $U_\al$, the fibre $\cH_x$
($x\in U_\al$) can be identified with $\bH$.
Using the local trivialisations of $P$, if $x\in U_\al\cap U_\beta$, the two
identifications of $\cH_x$ with $\bH$ from $U_\al$ and $U_\beta$ are related
by $\tilde\rho(\tilde g_{\al\beta}(x))\in U(\bH)$.
This shows the result that $\cH$ is associated to $Q$ and thus we have
the bundle isomorphism $\cH\cong Q\times_{LG}\bH$.

We now show that $\PP(\cF)$ defined in Sect.~\ref{sect:index} is an associated
bundle of $Q$ by $\hat\rho_0$.
We assume that on each $U_\al$, we can choose
$\lam_\al\not\in\text{spec}(D_x)$ for all $x\in U_\al$.
We then have a bundle of Fock spaces $\cF_{\lam_\al}$ over $U_\al$.
Since $\cH$ is an associated bundle of $Q$ by the representation
$\tilde\rho$ of $LG$ on the typical fibre $\bH$, $\PP(\cF_{\lam_\al})$
is an associated bundle of $Q_\al\to U_\al$ by the representation
$\hat\rho_{\lam_\al}\colon LG\to PU(\bF_{\lam_\al})$, where $\bF_{\lam_\al}$
is the Fock space of $\bH$ with a polarization corresponding to $\lam_\al$.
Under the natural identification of $\PP(\bF_{\lam_\al})$ and $\PP(\bF_0)$,
the representation $\hat\rho_{\lam_\al}$ is the same as $\hat\rho_0$.
This crucial step follows from the independence of representation
of $LG$ on the cuts, as explained in Sect~\ref{sect:Bog}.
Therefore, $\PP(\cF)$ is an associated bundle of $Q$; in fact it is
isomorphic to $\PP(\cF^\rho)$.

Since the line bundles $\cL\to\PP(\cF)^{[2]}$ in Sect.~\ref{sect:index} and
$\cL^Q\to\PP(\cF^\rho)^{[2]}$ in Sect.~\ref{sect:caloronb} are both constructed
from the universal line bundle, we have the following results.

\begin{thm}
Let $P\to S^1\times X$ be a principal $G$-bundle and let $Q\to X$ be the
corresponding principal $LG$-bundle.
Let $\rho$ be a finite dimensional unitary representation of $G$ on $V$ and
let $E=P\times_GV$.
Let $\hat\rho_0$ be the projective representation of $LG$ on $\bF_0$, the
Fock space of $\bH=L^2(S^1,V)$ with polarization at $0$.
Let $\cH$ be the Hilbert bundle whose fiber $\cH_x$ at $x\in X$ is the
$L^2$-completion of $\Gam(S^1\times\{x\},E)$.
Then the projectivized Fock bundle $\PP(\cF)$, constructed from $\cH$,
is isomorphic to the associated bundle
$Q\times_{\hat\rho_0}\PP(\bF_0)=\PP(\cF^\rho)$.
Furthermore, the index bundle gerbe of the family of Dirac operators on
$S^1$ coupled to $E$ is isomorphic to the caloron bundle gerbe.
\end{thm}

An alternative, less explicit argument for the same result goes as follows.
For simplicity, we will assume for the rest of the section that $H^3(X,\ZZ)$
is torsion free, i.e., the torsion subgroup of $H^3(X,\ZZ)$ is trivial.
It is well known that fibre bundles over $X$ whose typical fibre is
a projective (infinite dimensional) Hilbert space are classified up to
isomorphism by their Dixmier-Douady classes in $H^3(X,\ZZ)$.
Conversely, every class in $H^3(X,\ZZ)$ is the Dixmier-Douady class of a
fibre bundles over $X$ whose typical fibre is a projective (infinite
dimensional) Hilbert space.
This is essentially contained in Proposition~2.1 of \cite{AS}.

Given a unitary representation $\rho\colon G\to U(V)$, $E=P\times_GV$ is a
hermitian vector bundle over $S^1\times X$.
A connection $\tilde A$ on $P$ induces a hermitian connection $\nabla^E$ on
$E$, whose curvature is $F^E$.
For the family of self-adjoint Dirac operators on $S^1$ coupled to $E$, the
index bundle gerbe $(\PP(\cF), \cL)$ is constructed in Sect.~\ref{sect:index},
and its Dixmier-Douady class in $H^3(X,\ZZ)$ is represented by the
$3$-curvature form on $X$ \cite{CM, CMM, L}
$$  
\frac{1}{8\pi^2}\int_{S^1}\tr_V\,(F^E\wedge F^E),  
$$
where $F^E$ is the curvature of the induced connection $\nabla^E$ on $E$.
This is equal to $H_\rho=\iota_\rho\,H$ according to Proposition~\ref{hom}.
Therefore the Dixmier-Douady class of $\PP(\cF)$ is equal to that of the
bundle gerbe $(Q^\rho,\cL^\rho)$.

On the other hand, the homomorphism $LSU(V)\to U_{\rm res}(\bH)$ pulls back 
the basic central extension of $U_{\rm res}(\bH)$ to that of $LSU(V)$
(see \cite{PS}, Sect.~6.7), whereas the Shale-Stinespring embedding
$\sig_0\colon U_{\rm res}(\bH)\to PU(\bF_0)$ pulls back the basic
central extension $U(\bF_0)\to PU(\bF_0)$ of $PU(\bF_0)$ to that of
$U_{\rm res}(\bH)$. 
The latter can be argued from the literature as follows.
One has the commutative diagram
\begin{equation} \label{eqn:linebdles}
\xymatrix{U_{\rm res}(\bH) \ar[r] \ar[d]_{} & PU(\bF_0) \ar[d]^{} \\
Gr_{\rm res}(\bH) \ar[r]^{} & \PP(\bF_0),}    
\end{equation}
where $Gr_{\rm res}(\bH)$ denotes the restricted grassmannian, which is the
quotient of $U_{\rm res}(\bH)$ by a contractible subgroup (cf.\ page 115 in
\cite{PS}). 
The vertical arrows in equation \eqref{eqn:linebdles} are the projection maps,
which have degree one since the fibres are contractible.
Similarly, since the fibre of $PU(\bF_0)\to\PP(\bF_0)$ is also contractible,
the Dixmier-Douady class of the $\PP(\bF_0)$ bundle is equal to that of the
$PU(\bF_0)$ bundle.
The bottom horizontal arrow in \eqref{eqn:linebdles} is the Pl\"ucker
embedding.
The sentence just after equation (7.7.4) on page 116 in \cite{PS} says that
the pullback of the tautological line bundle on projective Fock space
$\PP(\bF_0)$ (under the Pl\"ucker embedding) is the determinant line bundle
on the restricted grassmannian $Gr_{\rm res}(\bH)$. 
Page 115 in \cite{PS} identifies the determinant line bundle on the restricted
grassmannian $Gr_{\rm res}(\bH)$ with the central extension of restricted
unitary group $U_{\rm res}(\bH)$.
Therefore the Dixmier-Douady class of the bundle gerbe $(Q^\rho,\cL^\rho)$
is equal to that of the induced $PU(\bF_0)$ bundle and that of the projective
Hilbert space bundle $\PP(\cF^\rho)$.

Thus, both the projective bundle $\PP(\cF)$ on $X$ associated to the index
bundle gerbe and the projective bundle $\PP(\cF^\rho)$ on $X$ associated to the
caloron bundle gerbe have the same Dixmier-Douady class in de Rham cohomology.
By our assumption on $H^3(X,\ZZ)$ and by the classification theorem for
projective Hilbert bundles over $X$, there is an isomorphism of projective
Hilbert bundles over $X$, $\PP(\cF^\rho)\cong \PP(\cF)$.
In particular, there is an induced isomorphism of fibred products
$\PP(\cF^\rho)^{[2]}\cong \PP(\cF)^{[2]}$.
Finally, we conclude that the caloron bundle gerbe $(\PP(\cF^\rho),\cL^\rho)$
and the index bundle gerbe $(\PP(\cF),\cL)$ are isomorphic.

\section{Bundle gerbes on moduli spaces associated to Riemann surfaces}
\label{sect:moduli}

In this section we apply the constructions of bundle gerbes in
Sect.~\ref{sect:index} and Sect.~\ref{sect:caloronb} to construct bundle
gerbes on moduli spaces associated to Riemann surfaces whose Dixmier-Douady
classes are generators of the third integral cohomology groups of these
moduli spaces.
We first review the construction of moduli spaces and universal bundles.

Let $G$ be a compact, connected, semisimple real Lie group whose centre is
$Z(G)$.
Let $G_\ad$ be the quotient $G/Z(G)$, which has a trivial centre.

Let $\Sig$ be a compact, connected Riemann surface of genus $g>1$, and let
$\tilde\Sig$ be its universal covering space, with a right action of the
fundamental group $\pi_1\Sig$.
There is a left $G$-action on the set of homomorphisms $\Hom(\pi_1\Sig,G)$
given by $G\ni h\colon\phi\mapsto\Ad_h\circ\phi$.
An element $\phi$ is irreducible if the isotropy subgroup is $Z(G)$, that is,
$\Ad_h\circ\phi=\phi$ for $h\in G$ implies $h\in Z(G)$.
Let $\Hom\irr(\pi_1\Sig,G)$ be the subset of such.
The quotient $\cM=\Hom\irr(\pi_1\Sig,G)/G$ is the moduli space of flat
$G$-connections on $\Sig$.
More generally, fix an element $z\in Z(G)$.
Let $\Sig'$ be the surface $\Sig$ with a puncture.
Its fundamental group $\pi_1\Sig'$ is the central extension of $\pi_1\Sig$
by $\ZZ$.
Let $\Hom_z(\pi_1\Sig',G)$ be the set of homomorphisms such that the
generator of $\ZZ$ is mapped to $z$.
Again there is a $G$-action on it and we have the subset
$\Hom\irr_z(\pi_1\Sig',G)$ of irreducible homomorphisms.
The quotient $\cM_z=\Hom\irr_z(\pi_1\Sig',G)/G$ is the moduli space of
flat $G$-connections on $\Sig'$ whose holonomy around the puncture is $z$.

For any $z\in Z(G)$, there is a left $G$-action on
$\tilde\Sig\times\Hom\irr_z(\pi_1\Sig',G)\times G$ given by
$G\ni h\colon(p,\phi,g)\mapsto(p,\Ad_h\circ\phi,hg)$.
There is also a left $\pi_1\Sig$-action on the same space given by
$\pi_1\Sig\ni\gam\colon(p,\phi,g)\mapsto(p\gam^{-1},\phi,\phi(\gam)g)$.
It is easy to check that the two actions commute and hence there is a left
action of $G\times\pi_1\Sig$ on
$\tilde\Sig\times\Hom\irr_z(\pi_1\Sig',G)\times G$.
Let $\cU=(\tilde\Sig\times\Hom\irr_z(\pi_1\Sig',G)\times G)/G\times\pi_1\Sig$
be the quotient space and let $\pi\colon\cU\to(\tilde\Sig\times
\Hom\irr_z(\pi_1\Sig',G))/G\times\pi_1\Sig=\Sig\times\cM$ be the projection
given by $\pi\colon[(p,\phi,g)]\mapsto[(p,\phi)]=([p],[\phi])$, where
$\phi\in\Hom\irr_z(\pi_1\Sig',G)$, $p\in\tilde\Sig$ and $g\in G$.
Then $\cU$ is a principal $G_\ad$-bundle over $\Sig\times\cM_z$.
To see this, we note that there is a right $G$-action on each fibre
$\pi^{-1}([p],[\phi])=\set{[(p,\phi,g)]}{g\in G}$ given by
$G\ni h\colon[(p,\phi,g)]\mapsto[(p,\phi,gh)]$.
The isotropic subgroup is $Z(G)$ because if
$[(p,\phi,g)]=[(p,\phi,gh)]=[(p,\Ad_{ghg^{-1}}\circ\phi,g)]$, then
$ghg^{-1}\in Z(G)$ or $h\in Z(G)$ by the irreducibility of $\phi$.
The right $G$-action thus descends to a free action of $G_\ad$ on $\cU$.

The bundle $\pi\colon\cU\to\Sig\times\cM_z$ is called the universal
bundle~\cite{Jef} because for each $[\phi]\in\cM_z$, the restriction
of $\cU$ to $\Sig\times[\phi]$ is the flat $G_\ad$-bundle over $\Sig$
defined by the homomorphism $\al\circ\phi\in\Hom\irr(\pi_1\Sig,G_\ad)$,
where $\al\colon G\to G_\ad$ is the quotient map.
In fact, there is a $G_\ad$-equivariant 1-1 map from $\cU|_{\Sig\times[\phi]}
=\pi^{-1}(\Sig\times[\phi])=\{[(p,\phi,g)]:p\in\tilde\Sig,g\in G\}$ to
$\tilde\Sig\times_{\al\circ\phi}G_\ad=\{[(p,g)]:p\in\tilde\Sig,g\in G_\ad\}$
given by $[(p,\phi,g)]\mapsto[(p,\al(g))]$.
Moreover, for any $\sig\in\Sig$, the restriction of $\cU$ to $\sig\times\cM_z$
is isomorphic to the $G_\ad$-bundle $\Hom\irr_z(\pi_1\Sig',G)\to\cM_z$.
To see this, we note that $\cU|_{\sig\times\cM_z}=\pi^{-1}(\sig\times\cM_z)=
\set{[(p_0,\phi,g)]}{\phi\in\Hom\irr_z(\pi_1\Sig',G),g\in G}$ for a fixed
$p_0\in\tilde\Sig$ such that $[p_0]=\sig$.
Under the bijection
$[(p_0,\phi,g)]=[(p_0,\Ad_g^{-1}\circ\phi,1)]\mapsto\Ad_g^{-1}\circ\phi$,
the right $G$-action
$G\ni h\colon[(p_0,\phi,1)]\mapsto[(p_0,\phi,h)]=[(p_0,\Ad_h^{-1}\circ\phi,1)]$
on $\cU|_{\sig\times\cM_z}$ corresponds to
$G\ni h\colon\phi\mapsto\Ad_h^{-1}\circ\phi$ on $\Hom\irr_z(\pi_1\Sig',G)$.

We assume that the moduli space $\cM_z$ is closed.
When $G=SU(n)$, $\cM_z$ is closed if and only if $z$ generates the center
$Z(SU(n))\cong\ZZ_n$.
Then the first Pontryagin class $x$ of $\Ad\cU\to\Sig\times\cM_z$ is a
generator of $H^4(\Sig\times\cM_z,\QQ)$.
For $G=SU(n)$, $H^2(\cM,\QQ)$ is generated by the slant product
$[\Sig]\backslash x$ while $H^3(\cM_z,\QQ)$ is generated by
$[\gam]\backslash x$, where $\gam\colon S^1\to\Sig$ are smooth loops whose
homology classes generate $H_1(\Sig,\ZZ) \cong \ZZ^{2g}$ when $g>1$ \cite{EK}.
There is no torsion for $H^3(\cM_z,\ZZ)$ if $z$ generates $Z(SU(n))$
\cite{BBGN}.
Given a Riemannian metric on $\Sig$, the universal bundle has a canonical
connection~\cite{ASin} whose curvature is denoted by
$F^\cU\in\Om^2(\Sig\times\cM_z,\Ad\cU)$.
In de Rham cohomology, $x$ is represented by
$-\frac{1}{8\pi^2}\bra F^\cU,F^\cU\ket$.
Therefore, the generators $[\Sig]\backslash x$ and $[\gam]\backslash x$
are represented by
$$ -\frac{1}{8\pi^2}\int_\Sig\bra F^\cU,F^\cU\ket\quad\text{and}\quad
     H_{\rho,\gam}
   =-\frac{1}{8\pi^2}\int_{S^1}(\gam\times\id)^*\bra F^\cU,F^\cU\ket, $$
respectively.

When $G=SU(n)$ and $z$ generates $Z(G)$, we apply the constructions of the
index and caloron bundle gerbes to obtain bundle gerbes on $\cM_z$ whose
Dixmier-Douady invariants are the generators of the third de Rham cohomology
group.
Given a smooth loop $\gam\colon S^1\to\Sig$, let $P_\gam=(\gam\times\id)^*\cU$
be the principal $G_\ad$-bundle over $S^1\times\cM_z$.
With a finite dimensional unitary representation $\rho$ of $G_\ad$ on $V$, we
have a principal $SU(V)$-bundle $P_\gam^\rho\to S^1\times\cM_z$ and hence an
$LSU(V)$-bundle $Q_\gam^\rho\to\cM_z$ by caloron correspondence, with a line
bundle $\cL_\gam^\rho$ over $(Q_\gam^\rho)^{[2]}$.
The universal connection on $\cU$ pulls back to $P_\gam$ and induces one on
$P_\gam^\rho$.
The latter determines a connection and a Higgs field for $Q_\gam^\rho$.
In particular, the Dixmier-Douady class of bundle gerbe
$(Q_\gam^\rho,\cL_\gam^\rho)$ is represented by $\iota_\rho\,H_{\rho,\gam}$,
which is one of the $2g$ generators of $H^3(\cM,\RR)$ if $g>1$.

Moreover, there are two ways to construct a projective Fock space bundles
over $\cM_z$ of the same Dixmier-Douady class as above.
First, by considering the the family of self-adjoint Dirac operators on
$S^1$ coupled to the vector bundle $E_\gam=P_\gam\times_GV$, we get the
index bundle gerbe $(\PP(\cF_\gam), \cL_\gam)$ from Sect.~\ref{sect:index}.
Second, using the $LSU(V)$ bundle $Q_\gam^\rho$, we have the caloron bundle
gerbe $(\PP(\cF_\gam^\rho), \cL_\gam^\rho)$ from Sect.~\ref{sect:caloronb}.
They are isomorphic by Sect.~\ref{theisom} and their Dixmier-Douady class is
equal to that of $(Q_\gam^\rho,\cL_\gam^\rho)$.

The construction of natural bundle gerbes also applies to other moduli
spaces associated to Riemann surfaces, such as the Hitchin moduli space
and the monopole moduli space.

\section{Conclusions and outlook}

In this paper, we have constructed bundle gerbes on moduli spaces $\cM$ of
flat $G$-bundles on a compact Riemann surface $\Sig$ of genus $g>1$, where
$G$ is a compact connected simply-connected Lie group. 
These are the index bundle gerbe and the caloron bundle gerbe, which we show
are isomorphic (not just stably isomorphic).
We have constructed a bundle gerbe connection and curving (or B-field), and
have computed the 3-curvature which represents the Dixmier-Douady class of
the bundle gerbe in de Rham cohomology.
The construction is such that it extends without change to other moduli spaces,
such as the moduli space of principal $U(n)$-bundles with fixed determinant
bundle and the moduli space of Higgs bundles.

In Sect.~\ref{sect:index}, given the basic setup, we have constructed the
index bundle gerbe, refining a construction of Segal.
It remains to define a natural bundle gerbe connection and curving on it,
and to compute the $3$-curvature form.
This has been established in the paper in a special case using 
Sect.~\ref{theisom}, which is the case that applies to the moduli spaces
considered here.

In the setting of Sect.~\ref{sect:caloron}, let $Q$ be a principal $LG$-bundle
over $X$ and $\hat\rho_0\colon LG\to PU(\bF_0)$ a positive energy
representation of $LG$. 
Then we can form the algebra bundle 
$$  \cK_Q=Q\times_{LG}\cK(\bF_0)\longrightarrow X  $$
with fibre the algebra of compact operators on Fock space, $\cK(\bF_0)$.
The space of continuous sections of $\cK_{Q,\hat\rho_0}$ vanishing at infinity,
$C_0(X, \cK_{Q,\hat\rho_0})$, is a $C^*$-algebra, which is a continuous trace
algebra with spectrum equal to $X$.
This non-commutative algebra is only locally Morita equivalent to continuous
functions on $X$.
The operator $K$-theory $K_\bullet(C_0(X,\cK_{Q,\hat\rho_0}))$ is the twisted
K-theory $K^\bullet(X,H)$ (cf.~\cite{Ros}), where $H$ is the Dixmier-Douady
class of $\PP(\cF^\rho)$ (see \cite{BCMMS} for a alternate description of
twisted $K$-theory). 
In particular, for each of the $2g$ bundle gerbes that we have constructed on
the moduli space $\cM$, we can form the twisted $K$-theory, and it would be
interesting to compute these.

Another interesting problem is to give natural constructions of holomorphic
bundle gerbes on the moduli space stable holomorphic bundles over a Riemann
surface $\Sig$.
If the degree and the rank are coprime, then the $2g$ generators of the third
cohomology are of types $(2,1)$ and $(1,2)$ for $g>1$~\cite{EK2}.
According to section~7 in \cite{MS03}, there are hermitian holomorphic bundle
gerbes that have these as Dixmier-Douady classes, but the challenge is to find
natural constructions for them.

Next we outline a construction of bundle gerbes on other moduli spaces
such as the moduli space of anti-self-dual (ASD) connections on a compact four
dimensional Riemannian manifold $M$ such that $\dim(H_1(M,\RR))>0$. 
More precisely, let $P$ denote a principal $G$-bundle over $M$, where $G$ is
a compact semisimple Lie group, and let $ \cA_P^-$ denote the space of all ASD
connections on $P$.
We will assume that $\cA_P^-\ne\emptyset$, which occurs under various
hypotheses on $P$ and $M$, cf. the introduction in~\cite{Taubes}.

If $\cG_P$ denotes the gauge group of $P$, then one has the universal
principal $G_\ad$ bundle
$$
G_\ad\longrightarrow(P\times\cA_P^-)/\cG_P\longrightarrow M\times\cA_P^-/\cG_P.
$$
Using this, we may also construct our bundle gerbes, for instance as follows.
For $\gam\colon S^1\to M$ a generator of $H_1(M,\ZZ)$, consider the restricted
bundle
\begin{equation}\label{princ}
G_\ad\longrightarrow(\gam^*P\times\cA_P^-)/\cG_P\longrightarrow
S^1\times\cA_P^-/\cG_P,
\end{equation}
which as before determines the caloron bundle gerbe over the ASD moduli space.
$$
LG_\ad\longrightarrow Q_\gam\longrightarrow\cA_P^-/\cG_P.
$$
One can similarly define the index bundle gerbe and show that these bundle
gerbes are isomorphic to each other and that their Dixmier-Douady is class
given by
$$      \int_{S^1}p_1\left(\gam^*(P)\times\cA_P^-)/\cG_P\right)
        \in H^3(\cA_P^-/\cG_P,\ZZ).              $$
Here $p_1\left(\gam^*P\times\cA_P^-)/\cG_P\right)$ denotes the first Pontryagin
class of the principal bundle in equation \eqref{princ}.


\bibliographystyle{abbrv}

\end{document}